\documentclass{eptcs}
\providecommand{\event}{ADG 2021} 

\usepackage{url}
\usepackage{xurl}
\usepackage{cite}
\usepackage{graphicx}
\usepackage{framed}
\usepackage{amsfonts}
\usepackage{listings}
\usepackage{color}
\usepackage{subfigure}
\usepackage{multicol}
\usepackage{amsmath}
\usepackage{breqn}

\usepackage{xcolor} 

\def\orcidID#1{\unskip$^{\mbox{\href{https://orcid.org/#1}{\scriptsize{[#1]}} }}$}

\def\titlerunning{Supporting Proving and Discovering Geometric Inequalities in GeoGebra by using Tarski}
\def\authorrunning{Brown, Kov\'acs \& Vajda}

\begin{document}

\title{Supporting Proving and Discovering Geometric Inequalities in GeoGebra by using Tarski}

\author{Christopher W.~Brown\orcidID{0000-0001-8334-0980}
\institute{United States Naval Academy\\
Annapolis, MD, USA}
\email{wcbrown@usna.edu}
\and
Zolt\'an Kov\'acs\orcidID{0000-0003-2512-5793}
\institute{The Private University College of \\
Education of the Diocese of Linz \\
Linz, Austria}
\email{zoltan@geogebra.org}
\and 
R\'obert Vajda\orcidID{0000-0002-2439-6949}
\institute{Bolyai Institute, University of Szeged \\
Szeged, Hungary}
\email{vajda@math.u-szeged.hu}
}

\maketitle              

\begin{abstract}
We introduce a system of software tools that can automatically prove
or discover geometric inequalities. The system, called \emph{Geo\-Gebra Discovery},
consisting of an extended version of \emph{Geo\-Gebra}, a controller web service
\emph{realgeom}, and the computational tool \emph{Tarski} (with the
extensive help of the \emph{QEPCAD B} system) successfully solves several
non-trivial problems in Euclidean planar geometry related to inequalities.


\end{abstract}

\section{Introduction}

For geometric \textit{equalities} there are already integrated tools in some dynamic geometry systems (DGS).
A well-known application is GeoGebra---it provides an automated reasoning toolset (ART)
that supports symbolic check of equality of
lengths of line segments in a geometric configuration, or, the equality of expressions.
Both features are supported by the \textit{Relation} tool or command,
or the low-level commands \textit{Prove} and \textit{ProveDetails} directly \cite{RelTool-ADG2014}.
GeoGebra's ART successfully proves hundreds of well-known theorems in
elementary planar geometry (see \url{https://prover-test.geogebra.org/job/GeoGebra_Discovery-provertest/66/artifact/fork/geogebra/test/scripts/benchmark/prover/html/all.html}
for a recent benchmark of 291 test cases), however, all of these results are related
to equalities.

In this extended abstract we describe our current work that focuses on \textit{inequalities}.
In our work we partly use the same algebro-geometric theory that
plays a crucial role in GeoGebra's toolset, based on the revolutionary work of Wu \cite{Wu78}, Chou \cite{Chou88},
and improved later by Recio and V\'elez \cite{RecioVelez99} with elimination theory.
This is the required algebraic basis for gaining a conjecture on \textit{a possibly fixed ratio of two quantities},
if equality does not hold between them. For example, the centroid of a triangle divides the medians in
ratio $2:1$, so here there is no equality between the parts, but a fixed a ratio, namely 2.
Once the elimination method is able to provide the exact ratio, the obtained conjecture becomes a mathematically proven proposition.

On the other hand, we partly use general purpose real quantifier elimination (RQE)
methods to find the best possible geometric constants in the related inequalities
if two expressions do not have a fixed ratio. Our implementation uses cylindrical algebraic decomposition (CAD)
that promotes effective RQE. As a special case, we can also directly prove geometric inequalities.
For example, theorems like the \textit{triangle inequality} or the inequality $s\leq3\sqrt3 R$
(where $s$ stands for the semiperimeter and $R$ for the circumradius) can already be mechanically
proven by our toolset in an intuitive way.

In this communication we do not go into the further detail on the hidden technical difficulties, but refer
to two recent papers: \cite{mcs2020} explains
the mathematical background on obtaining ratios in regular polygons,
while \cite{scsc-2020} focuses on the RQE related issues,
and points to a large set of benchmarks based on our tool. Instead, we focus
on the practical use: how our work can be fruitful for the student and the teacher in a classroom.
These are the two possible ways how
\textit{GeoGebra Discovery} obtains proofs of geometric inequalities:
\begin{itemize}
\item In the case of exploration it tries to solve the problem with elimination first,
and if this step is unsuccessful, then the given RQE problem setting will be outsourced
to the external tool \textit{realgeom}. We highlight that realgeom can work together with Wolfram's \textit{Mathematica}
(in this case the end user is expected to have a Mathematica subscription or to own
a \textit{Raspberry Pi} system that offers free access to Wolfram's tool, or use \texttt{wolframscript} that
offers a light version of the Mathematica kernel), or the free systems \textit{QEPCAD B} and
\textit{Tarski}.
\item In the case of request for a direct proof of an inequality, at the moment, the Tarski system is used.
\end{itemize}

Our experimental system is already capable
of solving a large set of open questions in planar Euclidean geometry.
Even if the mathematical methods we use are mostly well-known,
a full implementation of the various techniques in a single graphical application---being
freely available for millions of potential users---is completely new.

As precursors of our contribution we refer the reader to some related work. These include
the software packages SCRATCHPAD II \cite{scratchpad2} (the first uniform implementation
of the geometry prover prototypes) and
Java Geometry Explorer \cite{Ye_2011} (a program with graphical user interface
that can give animated explanations of proofs of \textit{equation based} geometry statements),
and the BOTTEMA solver \cite{bottema-program} (which is able to prove several inequalities
from \cite{Bottema69}); and the paper \cite{jianguo-zhang} on a survey of human readable proofs.
Further precursors of our work can be found in \cite{YangHouXia99,ChouGaoZhang2000,Geother,Xia2007,raglib}.

\section{Multiple ways for users to enter input}\label{multi}

School curriculum usually includes several relationships between expressions in a planar construction.
Every student has to learn the Pythagorean theorem which is an equality of two expressions,
namely $a^2+b^2$ and $c^2$ where $a$, $b$ and $c$ correspond to the lengths of sides in a right triangle.
In general, however, a simple question is not discussed: what happens if we omit the assumption
that the triangle has a right angle? By using our tool it can be quickly shown mechanically (here ``mechanically''
means that a \textit{machine} works in the background) that in general the
equality $\displaystyle{a^2+b^2>\frac{c^2}2}$
holds.

GeoGebra Discovery offers several ways to obtain this result. All of them require the first step that
a triangle $ABC$ must be drawn with GeoGebra's toolset.
GeoGebra automatically labels the sides $a$, $b$ and $c$.
After this, the user has the following options:
\begin{enumerate}
\item Directly proving the inequality by typing
  \texttt{Prove($a^2+b^2>c^2/2$)}. This is simple to input, but of
  course the user would have to know or guess that this inequality is
  the right one to try.

It is possible to arrive at the inequality by through multiple attempts, e.g.~\texttt{Prove($a^2+b^2>c^2$)}
(which is false), or \texttt{Prove($a^2+b^2>c^2/3$)} (which is true),
and then manually finding the largest
value $\mu$ such that $a^2+b^2>\mu\cdot c^2$ holds in general.

\item By entering the command \texttt{Compare($a^2+b^2$,$c^2$)}
the computer is asked to search for such a $\mu$ value, or, 
more generally, to find sharp non-negative constants $m$ and $M$ such that
\begin{equation}
m\cdot(a^2+b^2)\underset{(=)}{<}c^2\underset{(=)}{<}M\cdot(a^2+b^2)
\label{p1}
\end{equation}
always holds.
This command immediately prints the result in both the Algebra and Graphics Views
of GeoGebra, in a raw textual form.

\item For most users there is a high-level method available.
  After typing \texttt{Relation($a^2+b^2$,$c^2$)}
a popup window is shown in the program by reporting that $a^2+b^2$ and $c^2$ numerically
differ in the particular
case for the given triangle. After then, the user needs to click ``More$\ldots$'' to obtain the
exact symbolic result that considers all possible positions for points $A$, $B$ and $C$.

Alternatively, the user may define the sum $a^2+b^2$ by typing it: GeoGebra will automatically denote it by $d$. Now, by using
the Relation tool
\raisebox{-.15\height}{\includegraphics[width=0.4cm]{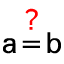}}
the user can directly compare $c$ and $d$. The underlying system learns that $c$ is a linear quantity
(in terms of some appearing lengths in the figure),
while $d$ is a quadratic one (again, in terms of lengths), so the investigated inequality should be of the form
\begin{equation}
m\cdot d^1\underset{(=)}{<}c^2\underset{(=)}{<}M\cdot d^1.
\label{p2}
\end{equation}
\end{enumerate}

These methods are supported in GeoGebra Discovery's latest versions \cite{geogebra-discovery}, but in the mainstream version not yet.

\section{Hierarchical structure of the system of tools}

\begin{figure}
\begin{center}
\includegraphics[scale=0.6]{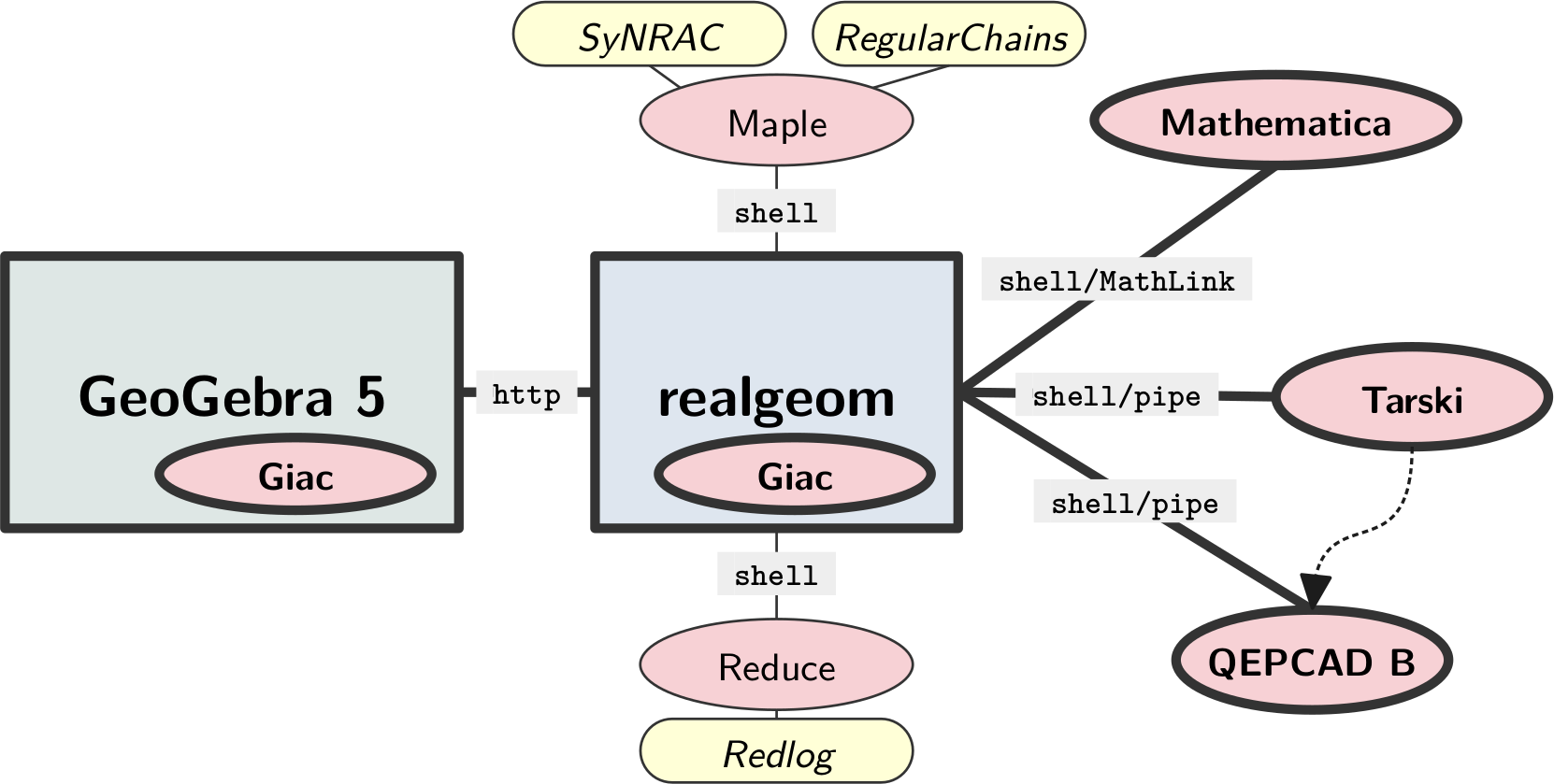}
\caption{Computation hierarchy for GeoGebra 5}\label{gg-rg3} 
\end{center}
\end{figure}

Fig.~\ref{gg-rg3} gives an overview of the underlying technical hierarchy. Our system
is modularized to be able to introduce other systems as optional parts at a later point.
Here we emphasize that it seems very difficult to start writing a completely new
computation system from scratch that provides all the mathematical background needed
to achieve our aims in a feasible timelimit for our input problems.
That is, re-using existing systems seems to be a much more reasonable way.

The
boldface items are used in reality to obtain the sharp constants $m$ and $M$, described
in the formulas above. That is, the current technical hierarchy
allows using \textit{Maple/SyNRAC} \cite{synrac2014}, Maple/\textit{RegularChains} \cite{ChenMaza16},
 and \textit{Reduce/Redlog} \cite{DolzmannSturm97} in the near future,
but for the moment only Mathematica's \cite{Mathematica} RQE subsystem,
\textit{QEPCAD B} \cite{CollinsHong91,Brown03} and Tarski  \cite{ValeEnriquez-Brown}
are available for the end user.
GeoGebra internally uses the Giac computer algebra system to perform symbolic computations
\cite{GiacGG-RICAM2013}.

Here we highlight that an important amount of work was performed to make our system work
in educational applications. First of all, each discussed component is free software. Second, we ported the QEPCAD B
and Tarski systems to Windows and made them available for the most recent macOS versions
to reach a large amount of users including teachers and students. The system is available
at \cite{geogebra-discovery} for all three majors platforms
(for Linux a 64 bit version and a Raspbian variant are published).

\section{A detailed example}

Let us consider a regular pentagon. We are interested in how far the points of two sides are
from each other.

\begin{figure}
\centering
\begin{multicols}{2}
 \includegraphics[scale=0.6]{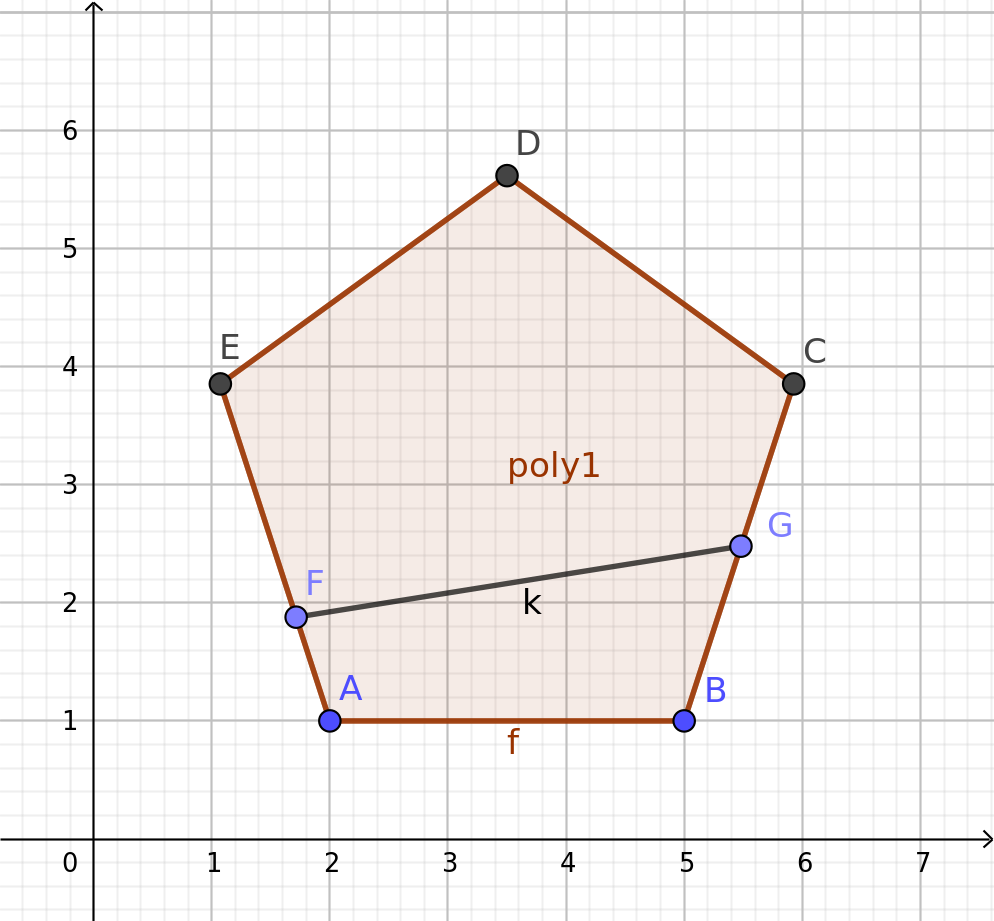}\vfill\null\columnbreak
 \begin{center}\includegraphics[scale=0.25]{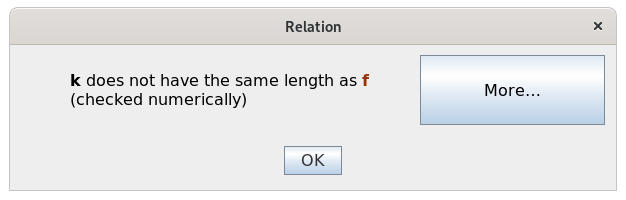}\end{center}
 \begin{center}\includegraphics[scale=0.25]{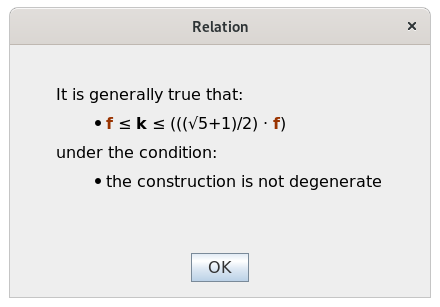}\end{center}
\end{multicols}
\caption{An example with a regular pentagon}\label{r5} 
\end{figure}

To formulate this precisely, we consider Fig.~\ref{r5} (left). A regular pentagon $ABCDE$ is given,
and on side $AE$ and $BC$, respectively, two points: $F$ and $G$. We connect them by a segment $k$.
We denote the segment $AB$ by $f$. Question: Is there a general relationship between $f$ and $k$?

When focusing on inequalities, this question can be split into two questions, namely:
\begin{enumerate}
\item What is $m$ such that $m \leq k/f$, and
\item what is $M$ such that $k/f \leq M$?
\end{enumerate}
The first question seems to be simple to answer: $m=1$, and this case can be reached when $F=A$ and $B=G$.
To find $M$ in the second question, however, some deeper considerations may be needed.

Here we emphasize that points $A$ and $B$ are freely chosen, and $F$ and $G$ have a one-dimensional freedom.
For the user, the easiest way to get a symbolic answer is the method we described in Section \ref{multi}, item 3:
by clicking on the Relation tool and then on each segment of $f$ and $k$ we run into a numerical
comparison first (Fig.~\ref{r5}, right, top). Then, by clicking on ``More$\ldots$'' we get, after a couple of seconds,
the result shown in Fig.~\ref{r5} (right, bottom).

\begin{figure}[ht!]
\begin{center}
\includegraphics[width=1\linewidth]{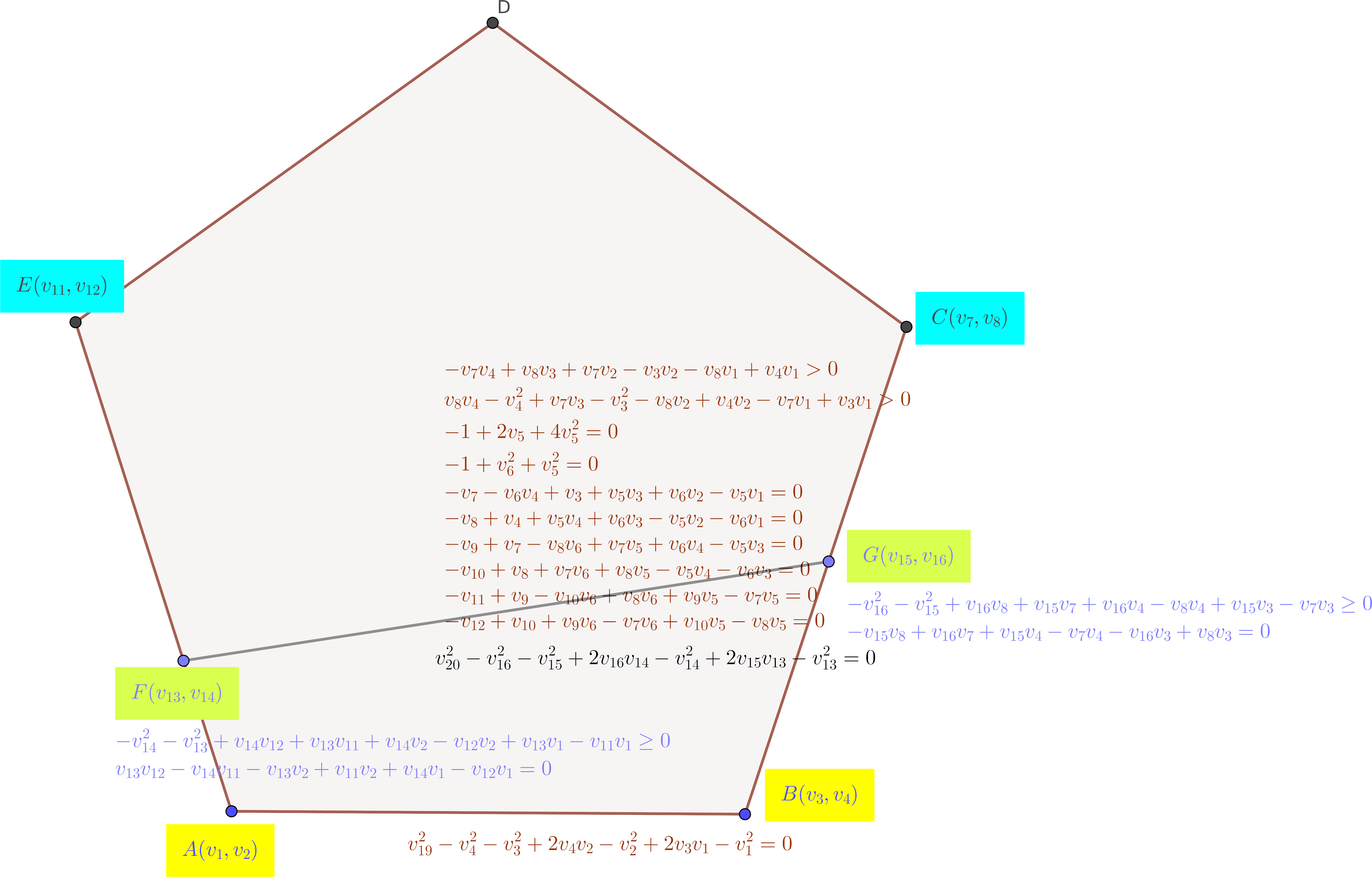}
\caption{Translation of the figure into a semi-algebraic system}\label{adg2-ex-expl} 
\end{center}
\end{figure}

The result is as expected, but we highlight here that the translation of the geometric setup into
a semi-algebraic system is not completely straightforward. Fig.~\ref{adg2-ex-expl} shows how
the geometric figure is translated into a set of equations and inequalities in an automated way.

The translation is as follows. First free points $A$ and $B$ are designated by coordinates $(v_1,v_2)$
and $(v_3,v_4)$, respectively. Coordinates of points $C(v_7,v_8)$, $D(v_9,v_{10})$ and $E(v_{11},v_{12})$ are uniquely determined by
the equations and inequalities shown inside the pentagon, presented in red, but repeated here in the text as well.
That is,
\begin{itemize}
\item[$\bullet$] $-v_7v_4+v_8v_3+v_7v_2-v_3v_2-v_8v_1+v_4v_1>0$ designates $C$ to be in the ``upper'' half-plane of line $AB$,
\item[$\bullet$] $v_8v_4-v_4^2+v_7v_3-v_3^2-v_8v_2+v_4v_2-v_7v_1+v_3v_1>0$ designates $C$ to be in the ``right'' half-plane of
a line going through $B$ and perpendicular to $AB$,
\item[$\bullet$] $-1+2v_{5}+4v_{5}^{2}=0$ defines the minimal polynomial of $\cos(2\pi/5)$ in the variable $v_5$, that is,
one of the roots of this polynomial is $\cos(2\pi/5)$,
\item[$\bullet$] $-1+v_{6}^{2}+v_{5}^{2}=0$ defines $v_6$ to express $\sin(2\pi/5)$.
\item[$\bullet$] The equations $-v_{7}-v_{6}v_{4}+v_{3}+v_{5}v_{3}+v_{6}v_{2}-v_{5}v_{1}=0$,
 $-v_{8}+v_{4}+v_{5}v_{4}+v_{6}v_{3}-v_{5}v_{2}-v_{6}v_{1}=0$, $-v_{9}+v_{7}-v_{8}v_{6}+v_{7}v_{5}+v_{6}v_{4}-v_{5}v_{3}=0$,
 $-v_{10}+v_{8}+v_{7}v_{6}+v_{8}v_{5}-v_{5}v_{4}-v_{6}v_{3}=0$, $-v_{11}+v_{9}-v_{10}v_{6}+v_{8}v_{6}+v_{9}v_{5}-v_{7}v_{5}=0$
and $-v_{12}+v_{10}+v_{9}v_{6}-v_{7}v_{6}+v_{10}v_{5}-v_{8}v_{5}=0$ designate the coordinates of the remaining points
of the regular pentagon. See \cite{mcs2020} for more details on this technique.
\end{itemize}
To put point $F(v_{13},v_{14})$ on line $AE$, we set the collinearity condition 
$v_{13}v_{12}-v_{14}v_{11}-v_{13}v_{2}+v_{11}v_{2}+v_{14}v_{1}-v_{12}v_{1}=0$,
and also the inequality
$-v_{14}^2-v_{13}^2+v_{14}v_{12}+v_{13}v_{11}+v_{14}v_2-v_{12}v_2+v_{13}v_1-v_{11}v_1\geq0$ is set to
ensure that $F$ lies between $A$ and $E$. An analogue equation and inequality is set up for
point $G(v_{15},v_{16})$.

In this example we skip variables $v_{17}$ and $v_{18}$ (they were used for an auxiliary computation),
but mention $v_{19}$ and $v_{20}$ that correspond to the length of $AB$ and $FG$, respectively.
They are defined by the polynomials
$v_{19}^2-v_4^2-v_3^2+2v_4v_2-v_2^2+2v_3v_1-v_1^2=0$ and
$v_{20}^2-v_{16}^2-v_{15}^2+2v_{16}v_{14}-v_{14}^2+2v_{15}v_{13}-v_{13}^2=0$, respectively.

Let us emphasize that we will focus on determining the possible values of $m$ (and $M$) by setting up
an equality among $v_{19}$, $v_{20}$ and $m$. This will be performed in step 2 (see below).

Having successfully translated the geometric figure in a set of semi-algebraic hypotheses $H=H_e\cup H_i$
(where $H_e$ stands for the set of equalities and $H_i$ for the inequalities), we consider the following steps:
\begin{enumerate}
\item We put $A$ and $B$ to $(0,0)$ and $(1,0)$, respectively, by learning that this can be done
without loss of generality.
\item It is possible that there are a finite number of candidate
  solutions that we can compute without considering $H_i$, so, roughly speaking,
  we construct the elimination ideal $(H_e\cup\{m\cdot v_{19}-v_{20}\})\cap\mathbb{Q}[m]$,
  and by using 
Gröbner basis computations, we hope for a zero-dimensional result. (Technically,
this step is done in GeoGebra, by using Giac.) Here this is not the case, so we go further.
\item To prepare the data for further processing, we \textit{delinearize} $H_e$. Roughly speaking,
this means that we try to eliminate as many variables as possible to prepare for the next step.
In our case we obtain the new equalities $H_e'$:
\begin{align}
-v_{13}^2+2 v_{13} v_{15}-v_{15}^2-v_{14}^2+2 v_{14} v_{16}-v_{16}^2+v_{20}^2&=0,\\
v_{10} v_7-2 v_7 v_8-v_{12}+v_8 v_9+v_8&=0,\\
-v_{10} v_8+v_8^2-v_{11}-v_7^2+v_7 v_9+v_7&=0,\\
v_7^2-2 v_7+v_8^2&=0,\\
4 v_7^2-6 v_7+1&=0,\\
-v_{15} v_8+v_8+v_{16} v_7-v_{16}&=0,\\
-v_{10}+2 v_7 v_8-v_8&=0,\\
v_7^2-v_7-v_8^2-v_9+1&=0,\\
-v_{11} v_{14}+v_{12} v_{13}&=0,\\
-m+v_{20}&=0,\\
-v_{19}+1&=0.
\end{align}
This step is performed in the realgeom tool.
\item By creating a quantified formula, here
\begin{equation}
\exists_{v_7,v_8,\ldots,v_{20}} (H_e' \land H_i \land v_{19}>0 \land v_{20}>0),
\end{equation}
we compute an equivalent, quantifier-free formula, which is
\begin{equation}
m > 0 \land m - 1 \geq 0 \land m^2 - m - 1 \leq 0.
\end{equation}
This step is computed with the help of the Tarski system. It reformulates the input
into a form that is hopefully easy to handle by QEPCAD B.

The quantifier elimination problems generated by this process are all
existentially closed except for the free variable $m$, and generally have
a number of equations in which some variables occur linearly.  These problems
are attacked via a script in Tarski's input language that does the following:
\begin{enumerate}
   \item applies quick simplification (command \texttt{bbwb});
   \item performs linear substitutions with case splitting to handle vanishing
      denominators and further quick simplification to eliminate as many variables
      as possible (command \texttt{qfr})---note that the resulting formula is a disjunction
      of existentially quantified conjunctions---and
   \item formulates each quantified conjunction for QEPCAD B and sends it
      to be solved independently (command \texttt{qepcad-qe}).  The final solution is then
      the disjunction of each of the sub-solutions which, recall, are formulas in
      the variable $m$ alone.
\end{enumerate}

\item We turn the quantifier-free formula into an interval that is simpler to interprete:
\begin{equation}
m\geq1 \land m\leq\frac{\sqrt5+1}{2}.
\end{equation}
This step is done in the realgeom tool by utilizing Giac.
\item Finally, GeoGebra reformulates this result as seen in Fig.~\ref{r5} (right, bottom).
\end{enumerate}
Here is how much time is spent in each step: The algebraic translation requires 50 ms,
the substitution and elimination steps (1--2) need 32 ms. For delinearization (3) 27 ms is needed.
Computation of the quantifier-free formula (4) requires 2216 ms, its rewrite to an interval (5)
needs 1 ms, and the final communication (6) requires another 1 ms. As a total, 2327 ms were needed
on a recent workstation (Ubuntu Linux 20.04, 8$\times$Intel(R) Core(TM) i7-8550U CPU @ 1.80GHz, 32 GB RAM).


This is one example among about two-hundred working test cases that vary between
simple constructions and some more complicated ones. We provide a benchmark of 131 examples,
available at
\url{https://prover-test.geogebra.org/job/GeoGebra_Discovery-comparetest/113/artifact/fork/geogebra/test/scripts/benchmark/compare/html/all.html}
that include several entries from the well-known Bottema database \cite{Bottema69}.
Also, direct proofs for 46 problems are benchmarked in the database 
\url{https://prover-test.geogebra.org/job/GeoGebra_Discovery-proverrgtest/33/artifact/fork/geogebra/test/scripts/benchmark/prover-rg/html/all.html}.
Both lists contain several examples that can be directly used in a classroom situation,
either to get a conjecture (with proof), or a confirmation/disproof, if the students already made a conjecture.

\section{Conclusion}

According to our benchmark, the Tarski system manages to solve 117 and 35 input problems within 30 secs timeout, respectively,
outperforming Mathematica with the number of successful tests (117 and 33 successful results, respectively).

By comparing the performance of Mathematica and Tarski we see that both systems
are very fast. In our benchmark they usually deliver the solution far below 1 second. This supports the idea
that Tarski can be used to obtain non-trivial results, by allowing
a large scale of users to study inequalities in a planar Euclidean geometry construction.
In classroom situations speed is very important: young learners are not expected to wait minutes,
hours or even days to obtain a result. Fast results are however not expected in a scientific environment.

\begin{figure}
\begin{center}
\includegraphics[width=0.5\linewidth]{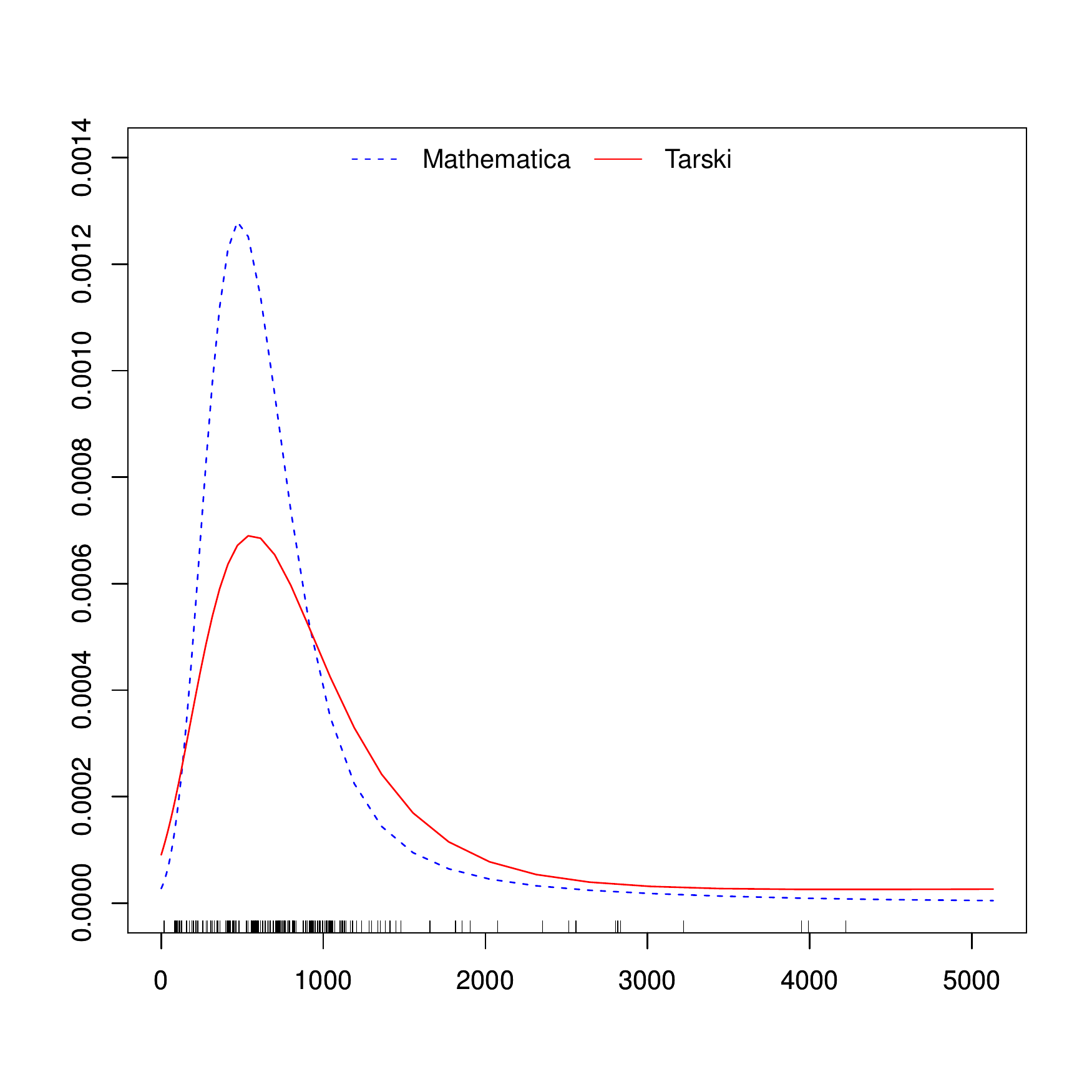}
\caption{Comparison of Mathematica (- - -) and Tarski (---)}\label{Mq} 
\end{center}
\end{figure}

We provide a figure that shows the density estimate of
time output of the first benchmark, visualized with statistics software R's
sm package \cite{BowmanAzzalini} on those 102 tests
that uniformly work on both systems (see Fig.~\ref{Mq}).
The $x$-axis shows the overall time in milliseconds. According to this diagram,
Mathematica's advantage is clear because of its smaller mean (1366 secs) and deviation (3485 secs),
however, Tarski's performance is also remarkable (2365 and 3940 secs, respectively).

In fact, more complicated inputs may result in getting no answer in feasible time in both systems.
It is well known that CAD has doubly exponential complexity in the number of variables
(see \cite{BrownDavenport2007, DavenportHeintz1988}), so in general there is no hope
to cover all problems in feasible time. But the questions we raise are quite specific,
so there seems to be much room for remarkable speedups.

We can consider the ``classical way'' of translating geometry into algebra as founded by
Wu and Chou in the 1980s, and its ``extension'' by adding non-equational constraints.
Surprisingly, just a couple of such constraints are required to faithfully describe
a planar Euclidean figure, if we consider the usual construction steps that are used
at secondary school level. These include
\begin{itemize}
\item the assumption of positive sign for the lengths (e.g.~those that belong to the sides of a triangle),
\item the \textit{betweenness} property if a point is attached to a segment (cf.~$F\in AE$ in Fig.~\ref{adg2-ex-expl}),
\item the condition to consider only the internal angle bisector (that is, to ignore the external ones),
\item as a special case of the previous one, the restriction that the center of the incircle of a triangle is a unique point (unlike 4 points that
include the centers of the excircles as well),
\item the constraint that a point is an inner point of a triangle, and
\item the uniqueness of a regular $n$-gon in the cases $n>2$, $n\neq4$ (unlike considering
some regular star-polygons at the same time, because they are indistinguishable in
the ``classical'' \textit{complex algebraic geometry setup}),
\end{itemize}
for the given hypotheses. These former inconveniencies can be handled properly in the ``extended'' \textit{real algebraic geometry setup},
and all of them are implemented in the most recent version of GeoGebra Discovery. 
In addition, we can avoid some difficulties that appeared when the thesis included a polynomial equality.
Formerly the thesis usually had to be rewritten to higher degree polynomials (e.g.~instead of $a+b=c$,
that is equivalent to $a+b-c=0$, the following equation had to be used: $(a+b-c)\cdot(b+c-a)\cdot(c+a-b)\cdot(a+b+c)=0$,
see \cite{MEP2018}). Real algebraic geometry successfully brings the possibility to forget these difficulties.

These extensions, however, have their price: higher computational complexity. To find an optimal balance
between ``imperfect but feasible'' (via complex algebraic geometry) and ``perfect but infeasible''
(via real algebraic geometry) seems to be a challenging research
for the next steps in our work.

\section{Future work}

Our final aim is to embed our experimental network of the four systems into one tool,
namely, into the \textit{mainstream version} of GeoGebra, developed at Johannes Kepler University of Linz,
Austria. It runs on several platforms, including desktop computers with various operating
systems (Windows, Mac and Linux), and the web platform, and also smartphones. This high
variety of technological background requires heavy simplification of the underlying technologies
and it calls for minimizing the applied external tools.

That is, for the long term our aim is to embed all external parts of our machinery into
GeoGebra. For the realgeom part, this seems to be achievable since it is a single Java
application with no graphical user interface: only some sophisticated algorithms
are included that can be copied to GeoGebra.

On the other hand, copying Tarski and QEPCAD B into
GeoGebra seems to be somewhat more difficult. QEPCAD is based on Saclib \cite{saclib-gh}, a C library that
offers CAS computations. QEPCAD and Tarski are written in C++. Thus, it seems a viable way to
compile all systems as a Java Native Interface (JNI) and include them as a dynamically
loadable object on all desktop platforms. This project is an on-going work.

Besides the native platform, a web application and a smartphone version should also be supported.
Compilation of C/C++ applications for the web is already well supported by various tools. Here
we highlight Emscripten \cite{emscripten} that provides multiple targets including JavaScript or WebAssembly.
Porting Tarski and QEPCAD B to the web is also an on-going work.

\section{Acknowledgements}
The second author was partially supported by a grant MTM2017-88796-P from the Spanish MINECO
(Ministerio de Economia y Competitividad) and the ERDF (European Regional Development Fund).
The third author was supported by the EU-funded Hungarian grant EFOP-3.6.1-16-2016-00008.

\bibliography{kovzol,external}

\end{document}